\documentclass[12pt]{article}
\usepackage{amstex, amsthm}
\newtheorem{theorem}{Theorem}[section]
\newtheorem{lemma}[theorem]{Lemma}
\newtheorem{proposition}[theorem]{Proposition}
\newtheorem{corollary}[theorem]{Corollary}

\newtheorem{definition}[theorem]{Definition}
\newtheorem{example}[theorem]{Example}

\newtheorem{remark}[theorem]{Remark}

\numberwithin{equation}{section}

\def \End{{\rm End}}

\def \deg {{\rm deg}}

\def \Aut{{\rm Aut}}

\def \Res{{\rm Res}}
\def \wt {{\rm wt}}

\def \Z{{\mathbb Z}}
\def \N{{\mathbb N}}
\def \C{{\mathbb C}}
\def \l{{\lambda}}
\def \h{{\frak h}}
\def \b{\beta}
\def \a{\alpha}
\def \bex{\begin{example}\label}
\def \eex{\end{example}}
\def \bl{\begin{lemma}\label}
\def \el{\end{lemma}}
\def \bc{\begin{corollary}\label}
\def \ec{\end{corollary}}
\def \bd{\begin{definition}\label}
\def \ed{\end{definition}}
\def \bt{\begin{theorem}\label}
\def \et{\end{theorem}}
\def \bp{\begin{proposition}\label}
\def \ep{\end{proposition}}
\def \br{\begin{remark}\label}
\def \er{\end{remark}}
\def \a{{\alpha}}
\def \pf{{\bf Proof. }}

\def \<{\langle}
\def \>{\rangle}

\def\1{{\bf 1}}
\def\zzz{{\mathbb Z}}

\def\hhh{{\hat {\frak h}}}
\def\fh{{{\frak h}}}

\def\res{{\hbox {Res}}}

\title{Representations of a class of lattice type vertex algebras}
\author{Stephen Berman{\thanks{
Supported by the Natural Science and Engineering Research Council of
Canada}}
\\{\small Department of Mathematics and Statistics, University of
Saskatchewan,}
\\{\small Saskatoon, Saskatchewan, Canada S7N 5E6}
\\
\\Chongying Dong {\thanks {Supported by NSF grant DMS-9987656 and
faculty
research funds granted by the University of California at Santa Cruz}}
\\{\small Department of Mathematics, University of California at Santa
Cruz, CA 95064}
\\
\\Shaobin Tan{\thanks{
Supported by the National Science Foundation of China}}
\\{\small Department of
Mathematics, Xiamen University, Xiamen 361005, China}}
\date{}
\begin{document}

\maketitle

\section { Introduction}
 
The lattice vertex algebras, [B1], [FLM], form one of the most 
important and fundamental classes of 
vertex algebras. Beginning with a root lattice of  simply laced type 
one
constructs the fundamental representation for the corresponding affine 
Kac-Moody  
Lie algebra and this turns out to be one of the most basic  examples of 
a lattice vertex algebra ([FK],[S]). The 
vertex algebra associated to the Leech lattice plays a fundamental role 
in the
 construction of the moonshine vertex operator algebra ([B1],[FLM]), 
while the vertex 
algebra associated to the rank 2 Lorentz lattice is used in 
constructing the Monster Lie algebra which in turn is used in a proof of the 
moonshine conjecture. Moreover lattice vertex algebras have been studied 
extensively from several points of view. Much work has been done on the 
representation 
theory ([FLM], [D], [DLM1], [DLM2]) and on fusion rules, [DL], for the 
lattice vertex 
algebras. Also, much is known about their automorphism groups, [DN], 
and there are some 
characterizations known [LX].
 
In this paper we study the representation theory for certain ``half 
lattice'' 
vertex algebras. Let $L$ be an even lattice and let $V_L$ be the 
associated lattice vertex algebra.  Then, as a vector space, $V_L$ is 
the 
tensor product of a symmetric algebra $S(\h \otimes_{\C}t^{-1} 
\C[t^{-1}])$ with a  
group algebra $\C[L]$ where $\h = \C \otimes_\Z L$.
 Thus,$V_L=S(\h \otimes_{\Z}t^{-1} \C[t^{-1}]) \otimes_{\C} \C [L]$. 
The lattice $L$ we consider in this 
paper is spanned by $c_i,d_i$ for $i=1,...,\nu$ with the 
$\Z$-bilinear form  determined  by
$(c_i,c_j)=(d_i,d_j)=0$ and $(c_i,d_j)=k \delta_{i,j}.$ Here $k$ is an 
arbitrary non-zero integer.
Our half lattice vertex algebra $V$ is then defined to be $S(\h 
\otimes_{\Z}t^{-1}\C[t^{-1}])\otimes_{\C}\C[L_C]$ where 
$L_C=\sum_{i=1}^{\nu}\Z c_i.$  
Notice that $V$ is not a lattice vertex algebra as originally defined 
in [B1] and [FLM], but it is a vertex subalgebra of $V_L$.

   The motivation for studying $V$, as defined above, comes from the 
representation  
theory of certain toroidal Lie algebras and certain other Lie algebras  
related to  
them. In [T] certain vertex operator representations were given for a 
Lie algebra 
which is constructed from a natural Jordan algebra using the TKK 
construction. This algebra is of great interest
in studying the structure theory and representation theory of extended 
affine Lie algebras in general   
because it has the smallest root system of any tame extended affine Lie 
algebra which is not of 
finite or affine type. It became clear that  to study the 
representation theory of 
this Lie algebra one should use techniques offered by taking a vertex 
algebra point of view and that a certain 
natural toroidal Lie algebra entered into the picture here. The 
representation theory of the toroidal algebras have been studied for some time 
now and recently they too have been viewed from the vertex algebra 
perspective (see [BBS] and the references therein). In fact, the work [BBS] 
studies and introduces a vertex algebra which is the tensor product of 
three other 
basic vertex algebras, and one of these is nothing but our $V$ above. 
The other two are VOA's associated to affine algebras and hence their 
representation theory is well developed. Thus, one needs the 
representation theory of the vertex algebra $V$ in order to understand that of the 
toroidal algebras in [BBS] as well as the algebras in [T]. This makes 
it natural to isolate $V$ and study it on its own.

Because of the motivation discussed above our main goal in this paper 
is to construct a large class of modules
for $V.$ Recall that in [FLM], for any even lattice $A$ 
and any element $\lambda$ in the dual lattice  of $A,$ 
an irreducible module $V_{A+\lambda}$ is constructed. The method in 
[FLM] for doing 
this uses the idea of coherent states and the authors make use of the 
fact that 
the lattice $L$ spans $\h$. In the work [D], employing the ideas of 
$Z$-algebras developed in [LW], 
it is proved that all irreducible $V_A$-modules are of the form $V_{A+ 
\lambda}$ as above. In our 
case, the sublattice $L_C$ does not span $\h$, and hence we cannot use 
the coherent state argument to 
 construct $V$-modules, but the $Z$-algebra techniques still play an 
important role.
 
The main tools we use to construct $V$-modules is the  theory of local 
systems, as
developed in [L2].
It turns out that our proofs also apply to the
``full lattice'' vertex algebra. Thus we have obtained another 
proof that $V_{A+\lambda}$ is a $V_A$-module.  More precisely,
we first construct an associative algebra $A$ and a large class of 
$A$-modules
in Section 3. We prove in Section 4 that for every $A$-module $W$ and
$\lambda\in\frac{1}{k}L_D$ where $L_D=\sum_{i=1}^{\nu}\Z d_i$ the
space $V_{\l,W}=S(\h \otimes_{\Z}t^{-1}\C[t^{-1}])\otimes_{\C}W$ is
a $V$-module. The definition of vertex operators $Y_{\l,W}(v,z)$
for $v\in V$ is essentially the same as in [FLM]. It follows from
the theory of local system [L2] that these operators 
form such a local space. It is well known that
the Jacobi identity for a module is equivalent to 
locality and associativity (see Section 4), so our main effort
is to prove these operators satisfy associativity. Unfortunately,
we cannot prove  this directly. However, as we show in Section 4,  
associativity follows from the 
fact that this local space is closed under all $n$th products 
$u(x)_nv(x)$ for 
$n\in \Z$ and $u,v\in V$ where $u(x)=Y_{\l,W}(u,x)$. Our results then 
follow from this. Finally, in Section 5,we also
discuss how to get $A$-modules from $V$-modules, and in Section 6 we 
study the 
Zhu algebra of $V$.

   All three authors take this opportunity to thank the Fields 
Institute in Toronto for their great hospitality while this work was being 
carried out.

\section{Vertex Algebra and Module}

In this paper, the sets of integers, positive integers and negative 
integers
will be denoted respectively by $\zzz, \zzz_+$ and $\zzz_-$. $z, z_1, 
z_2$
and $x$ will denote  formal variables. The elementary properties
of the $\delta$-function $\delta(z)=\sum_{n\in\zzz}z^n$ can be found in
[FLM] and [FHL]. First we give the definition
of a vertex algebra and  modules for them (see [B1], [FLM], [L1],
and we also recall the theory of local system developed in [L2].

\vspace{3mm}

\bd{dvoa}
A vertex algebra is a quadruple $(V,Y,{\bf 1},D)$ consisting of a
vector space $V$, a vector ${\bf 1}\in V,$ a linear map $Y$ from $V$ to
$({\rm End}V)[[z,z^{-1}]]$ and a linear map $D$ from
$V$ to $V$ satisfying the following axioms:

(1) For any $u,v\in V,$ $u_nv=0$ if $n$ is very large,

(2) $Y({\bf 1},z)=1$,

(3) $Y(u,z){\bf 1}\in V[[z]]$ and
$\lim_{z\rightarrow 0}Y(u,z){\bf 1}=u$ for $u\in V$,

(4) $[D, Y(u,z)]={\frac d {dz}}Y(u,z)=Y(Du,z)$ for $u\in V$,

(5) The Jacobi identity holds for $u,v\in V$:
\begin{eqnarray}
& 
&z_{0}^{-1}\delta\left(\frac{z_{1}-z_{2}}{z_{0}}\right)Y(u,z_{1})Y(v,z_{2})
-z_{0}^{-1}\delta\left(\frac{z_{2}-z_{1}}{-z_{0}}\right)Y(v,z_{2})Y(u,z_{1})
\nonumber\\
& 
&=z_{2}^{-1}\delta\left(\frac{z_{1}-z_{0}}{z_{2}}\right)Y(Y(u,z_{0})v,z_{2}).
\end{eqnarray}
\ed

We shall also use $V$ for the vertex algebra $(V,Y,{\bf 1},D)$.

\bd{dmodule}
Let $V$ be a vertex algebra.
A $V$-module is a vector space $W$ equipped with a  linear map
$Y_W$ from $V$ to
$({\rm End}\;W)[[z,z^{-1}]]$ satisfying the following axioms:

(1) For any $u\in V,$ $w\in W,$ $u_nw=0$ if $n$ is very large,

(2) $Y_W({\bf 1},z)=1,$

(3) The Jacobi identity holds for $u,v\in V$:
\begin{eqnarray}\label{2.2}
& &z_{0}^{-1}\delta\left(\frac{z_{1}-z_{2}}{z_{0}}\right)
Y_W(u,z_{1})Y_W(v,z_{2})
-z_{0}^{-1}\delta\left(\frac{z_{2}-z_{1}}{-z_{0}}\right)
Y_W(v,z_{2})Y_W(u,z_{1})
\nonumber\\
& &=z_{2}^{-1}\delta\left(\frac{z_{1}-z_{0}}{z_{2}}\right)
Y_W(Y(u,z_{0})v,z_{2}).
\end{eqnarray}
\ed

One of the important consequences of the definition of module is
the following $D$-derivative property:
$$Y_W(Du,z)=\frac{d}{dz}Y_W(u,z)$$
(see Lemma 2.2 of [DLM1]).

We now recall the theory of local system  of vertex operators from
[L2]. Let $M$ be a vector space. A {\it vertex operator} on $M$ is a 
formal series
$a(z)=\sum_{n\in {\Z}}a_{n}z^{-n-1}\in ({\rm End}\;M)[[z,z^{-1}]]$ such 
that
for any $u\in M$, $a_{n}u=0$ for sufficiently large $n$. All vertex
operators on $M$ form a vector space (over ${\C}$), denoted by $VO(M)$.
On $VO(M)$, we have a linear endomorphism $D={\frac{d}{dz}}$, the formal
differentiation.

Two vertex operators $a(z)$ and $b(z)$ on $M$ are said to be
{\it mutually local} if there is a
non-negative integer $k$ such that
\begin{eqnarray}
(z_{1}-z_{2})^{k}a(z_{1})b(z_{2})=(z_{1}-z_{2})^{k}b(z_{2})a(z_{1}).
\end{eqnarray}
A space $S$ of vertex operators is said to be {\it local} if any two
vertex operators of $S$ are mutually local, and a maximal local space 
of
vertex operators is called a {\it local system}.

Let $V$ be a local system on $M$. Then $V$ is closed under the formal
differentiation $D=\dfrac{d}{dx}$. For $a(x),b(x)\in VO(M)$, we define
\begin{eqnarray}\label{evp}
& &Y(a(x),z)b(x)\nonumber\\
&=&\Res_{z_{1}}\left(z^{-1}\delta\left(\frac{z_{1}-x}{z}\right)
a(z_{1})b(x)-z^{-1}\delta\left(\frac{x-z_{1}}{-z}\right)b(x)a(z_{1})\right)
\end{eqnarray}
Write $Y(a(x),z)=\sum_{n\in\Z}a(x)_nz^{-n-1}.$ Then
(\ref{evp}) is equivalent to
\begin{equation}\label{action}
a(x)_nb(x)=\Res_{z}((z-x)^na(z)b(x)-(-x+z)^nb(x)a(z))
\end{equation}
for $n\in \Z.$  Denote by $I(x)$ the identity endomorphism of $M$.

\bt{tli} [L2]
Let $M$ be a vector space and $V$ a local system on $M$. Then
$(V, Y, D, I(x))$ is a vertex algebra with $M$ as a natural module
such that $Y_M(a(x),z)=a(z)$ for $a(x)\in V.$
\et

Let $A$ be any local space of vertex operators on $M$. Then there 
exists a
local system $V$ that contains $A$. Let $\<A\>$ be the vertex 
subalgebra
of $V$ generated by $A$. Since the vertex operator ``product'' 
(\ref{evp})
does not depend on the choice of local system $V$, $\<A\>$ is 
canonical.
Then we have:

\bc{clocal} [L2] Let $M$ be a vector space and $A$ any local space of 
vertex
operators on $M$. Then $A$ generates a canonical vertex algebra $\<A\>$
with $M$ as a natural module such that $Y_M(a(x),z)=a(z)$ for
$a(x)\in A.$
\ec

Next we recall the well-known lattice vertex algebras from [B1]
and [FLM]. We are working in the setting of [FLM]. In particular, $L$ 
is a
lattice with nondegenerate symmetric ${\Bbb Z}$-bilinear ${\Bbb 
Z}$-valued
form $\langle\cdot,\cdot\rangle;$
 ${\frak h}=L\otimes_{\Bbb Z}{\Bbb C};$
$\hat{\frak h}_{\Bbb Z}$ is the corresponding Heisenberg algebra; 
$M(1)$
is the associated irreducible induced module for $\hat{\frak h}_{\Bbb 
Z}$ such that the canonical
central element of $\hat{\frak h}_{\Bbb Z}$ acts as 1; $(\hat L,{}^-)$ 
is
a central extension of $L$ by group $\langle
\kappa|\kappa^2=1\rangle $ with commutator map
$c(\alpha,\beta)=\kappa^{(\alpha,\beta)}$ for $\alpha,\beta\in L;$
$\chi$ is a faithful character of $\langle\kappa\rangle$ such that
$\chi(\kappa)=-1;$ ${\Bbb 
C}\{L\}=\mbox{Ind}_{\langle\kappa\rangle}^{\hat L}{\Bbb C}_{\chi}\simeq {\Bbb C}[L]$ (linearly), where ${\Bbb 
C}_{\chi}$ is the
one-dimensional $\langle\kappa\rangle$-module defined by $\chi;$
$\iota(a)=a\otimes 1\in {\Bbb C}\{L\}$ for $a\in \hat L;$ 
$V_L=M(1)\otimes {\Bbb C}\{L\};$ $\mbox{\bf 1}=\iota(1);$ 
$\omega=\frac{1}{2}\sum_{r\ge 1}\beta_r(-1)^2$
where $\{\beta_1,\beta_2,...\}$ is an orthonormal basis of ${\frak 
h};$
\begin{equation}\label{5.8}
\begin{array}{lcr}
V_L&\to& (\mbox{End}\,V_L)[[z,z^{-1}]]\hspace*{3.6 cm} \\
v&\mapsto& Y(v,z)=\displaystyle{\sum_{n\in{\Bbb Z}}v_nz^{-n-1}\ \ \ 
(v_n\in
\mbox{End}\,V_L)}
\end{array}
\end{equation}
where the vertex operator $Y(v,z)$ is defined in detail in [FLM]; the
space $V_L$ carries a natural ${\Bbb Z}$-grading determined
by the conditions wt$(1\otimes\iota(a))=\frac{1}{2}\langle\bar a,\bar
a\rangle.$ Then $(V_L,Y,{\bf 1},L(-1))$ is a simple vertex algebra
where $Y(\omega,z)=\sum_{n\in Z}L(n)z^{-n-2}$ (see [B1] and [FLM]).
Moreover, the component operators $L(n)$ satisfy the Virasoro
algebra relation with central charge rank$L$ and 
$V_L=\oplus_{n\in\Z}(V_L)_n$ is
$\Z$-graded where $(V_L)_n$ is the eigenspace of $L(0)$ with
eigenvalue $n$ (see [B1] and [FLM]).

In this paper we consider certain even lattices below and study
the modules for the ``half lattice'' vertex subalgebra of
$V_L.$

Let $L_C=\oplus^{\nu}_{i=1}\Z c_i$, $L_D=\oplus^{\nu}_{i=1}\Z d_i$
and $L=L_C+L_D$. Define a symmetric bilinear form $(\cdot)$ on $L$
such that
$$
(c_i,d_j)=k\delta_{ij},\;\;\;\;\;(c_i,c_j)=(d_i,d_j)=0
$$
for $1\le i,j\le\nu$, where $k\in \zzz\setminus\{0\}$ is a constant.
Then $L$ is an even lattice of rank $2\nu$, and $L$ is unimodular if 
$k=1.$

For $\alpha\in (L_C\otimes_{\Z}\C)/kL_C$ we define an automorphism
$g_{\alpha}\in \Aut V_L$ by
$$
g_{\alpha}(u\otimes\iota(a) )=e^{2\pi i(\alpha,\bar a)}u\otimes 
\iota(a)
$$
for $u\in M(1)$ and $a\in\hat L$ (see [DM1]). Then
$G=\{g_{\alpha}|\;\;\alpha\in  (L_C\otimes_{\Z}\C)/L_C\}$ is an
abelian subgroup of $\Aut V_L$. Let $V^G_L$ be the space
of $G$-invariants. Then
$$
V^G_L=M(1)\otimes\C\{L_C\}
$$
where $\C\{L_C\}$ is spanned by $\iota(a)$ for $a\in \hat L$ such that
$\bar a\in L_C.$ It is clear that the algebra $\C\{L_C\}$ is 
isomorphic
to the group algebra $\C[L_C].$ So we will use $e^{\alpha}$ for
 basis elements of  $\C[L_C]$ corresponding to $\alpha\in L_C.$

\bp{p2.1}
$V^G_L=\oplus^{\infty}_{n=0}V_n$ is a $\Z$-graded
simple vertex subalgebra of $V_L$, where
$$ (V_L^G)_n=V_L^G\cap (V_L)_n=M(1)_n\otimes \C[L_C]$$ 
where 
$$M(1)_n=\{ \a_1(-n_1)\cdots \a_s(-n_s)|\a_i\in{\frak h}, n_i>0,\sum^s_{i=1}n_i=n\}.
$$
\ep

\pf The gradation is clear from the definition of the grading of
$V_L.$ In order to see that $V_L^G$ is simple we note that
$$V_L=\oplus_{\beta\in L_D}V_L^{\beta}$$
where $V_L^{\beta}=\{v\in V_L|g_{\alpha}v=e^{2\pi i(\alpha,\beta)}v,
\alpha\in  (L_C\otimes_{\Z}\C)/kL_C.$ Clearly, $V_L^0=V_L^G$ and
$u_nV_L^{\beta}\subset V_L^{\beta+\gamma}$ for $\beta,\gamma\in L_D,$
$u\in V_L^{\gamma},$ $n\in\Z.$ Since $V_L$ is simple, $V_L$
is spanned by $u_nv$ for $u\in V_L$ and $n\in \Z$ where $v$ is any
nonzero vector in $V_L$ (see Corollary 4.2 of [DM2] and Proposition 4.1 
of [L1]). It is immediate now that $V_L^G$ is simple and
each $V_L^{\beta}$ is a simple $V_L^G$-module.
\qed

\section{Associative Algebra and Modules}

Our main goal in this paper is to study the representation theory for 
$V=V^G_L.$ It turns out that this is closely related to the 
representation  
theory of an associative algebra $A$ which we will define in this 
section. It 
is convenient for us to begin by studying a larger algebra $B$ which 
then has 
$A$ as a homomorphic image. The idea to consider these algebras comes 
from considering 
$Z$-algebras in our setting. Although our study was motivated by the 
idea of $Z$-algebras we 
 will see in Section 6 that $A$ is precisely the Zhu algebra $A(V)$, 
[Z]. The algebra $B$ is essentially a 
 twisted tensor product of a group algebra on the lattice $L_C$ and the 
tensor algebra 
of the space $\h_D=\C \otimes_{\Z}L_D$ where the elements $d_i$ act on 
$\C [L_C]$ as derivations.
We also let $\h_C= \C \otimes_{\Z} L_C.$
 
Formally, we let  $B$ be an associative algebra generated by 
$e_{\alpha}$ and $d_i$ for
$\alpha\in L_C$ and $1\le i\le \nu$, subject to the following relation
$$
e_0=1,\;\;\;\; e_{\alpha+\beta}=e_{\alpha}e_{\beta},\;\;\; 
d_ie_{\alpha}-
e_{\alpha}d_i=(d_i,\alpha)e_{\alpha}
$$
for $\alpha,\beta\in L_C$, $1\le i\le \nu$.

\vspace{3mm}
\bd{d2.5} A $B$-module $W$ is called a weight module if
$W=\oplus_{\lambda\in {\frak h}_C}W_{\lambda}$, where
$$
W_{\lambda}=\{ w\in W|d_iw=(\lambda, d_i)w, i=1,...,\nu\}.
$$
\ed

Regarding ${\frak h}$ as an abelian group, then the group algebra
$\C[{\frak h}]$ is a weight module for $B$ where $e_{\alpha}$ acts
by addition and $d_i$ acts on $e_{h}$ for $h\in {\frak h}_C$ as
$(d_i,h).$ It is easy to see that  $\C[L_C+\lambda]$ is a simple
module for any $\lambda\in {\frak h}_C.$

\bl{l2.6} Any weight $B$-module is semisimple. Moreover
all the simple $B$-modules (up to isomorphism) are  $\C[L_C+\lambda]$ 
for some $\lambda\in {\frak h}_C.$
\el

\pf Let $w\in W_{\lambda}$ be nonzero for some $\lambda\in {\frak 
h}_C$. Set
$\<w\>=\oplus_{\alpha\in L_C}\C e_{\alpha}w$. Then $\<w\>$ is a simple
$B$-submodule of $W$ isomorphic to $\C[L_C+\lambda]$ and the lemma 
follows.
\qed

\vspace{4mm}
Next we are going to construct a large class of simple B-modules which
contains the weight modules as a special case. For simplicity we assume
that $k=1$ (the unimodular case). Let $1\le\mu\le\nu+1$ be an
integer. Consider the Laurent polynomial algebras
$\C[t_1^{\pm 1},\cdots,t_{\mu-1}^{\pm 1}]$,
and $\C[t_{\mu},\cdots,t_{\nu}]$ with commuting variables. For any 
fixed
Laurent polynomials $f_i=f_i(t_1,\cdots,
t_{\mu-1})\in\C[t_1^{\pm 1},\cdots,t_{\mu-1}^{\pm 1}]$ for
$i=1,2,\cdots,\mu-1$, and nonzero complex numbers $a_i$ for
$i=\mu,\cdots,\nu$  let 
$\omega=\omega(f_1,\cdots,f_{\mu-1}|a_{\mu},\cdots
, a_{\nu})$ be a symbol. We will also use the symbols 
$\omega(|a_1,\cdots,
a_{\nu})$ if $\mu=1$, and $\omega(f_1,\cdots, f_{\nu}|)$ if $\mu=\nu+1$.
Define
$$
M_{\omega}=\C[t_1^{\pm
1},\cdots,t_{\mu-1}^{\pm 1},t_{\mu}\cdots,t_{\nu}]\omega.
$$
We define actions of $e_{\alpha}, d_j$ on $M_{\omega}$ for
$\alpha=\sum_{i=1}^{\nu}m_ic_i$ and $1\le j\le \nu$. For any $f\in
\C[t_1^{\pm 1},\cdots,t_{\mu-1}^{\pm 1},t_{\mu}\cdots,t_{\nu}]$,
$$
e_{\alpha}. 
f\omega=\left[(\prod^{\mu-1}_{i=1}t_i^{m_i})\prod_{i=\mu}^{\nu}
(a_ie^{-\partial_i})^{m_i}f\right]\omega,
$$
$$
d_j.
f\omega=(t_j\partial_jf+f_jf)\omega,\;\;\;\;{\hbox{for}}\;\;1\le j\le
\mu-1,
$$
$$
d_j. f\omega=t_jf\omega,\;\;\;\;{\hbox{for}}\;\;\mu\le j\le
\nu.
$$
where $\partial_i={\frac {\partial}{\partial t_i}}$.

\bt{t3.3} $M_{\omega}$ is a simple $B$-module.
\et
 
\pf We first prove that $M_{\omega}$ is an $B$-module. For $1\le
i\le\mu-1$, we have
$$
(d_je_{\alpha}-e_{\alpha}d_j).f\omega
=d_j.\left[(\prod^{\mu-1}_{i=1}t_i^{m_i})\prod_{i=\mu}^{\nu}
(a_ie^{-\partial_i})^{m_i}f\right]\omega-e_{\alpha}.(t_j\partial_j+f_j)f\omega
$$
$$
=(t_j\partial_j+f_j)\left[(\prod^{\mu-1}_{i=1}t_i^{m_i})\prod_{i=\mu}^{\nu}
(a_ie^{-\partial_i})^{m_i}f\right]\omega
$$
$$
-\left[(\prod^{\mu-1}_{i=1}t_i^{m_i})\prod_{i=\mu}^{\nu}
(a_ie^{-\partial_i})^{m_i}\right](t_j\partial_j+f_j)f\omega
$$
$$
=m_j\left[(\prod^{\mu-1}_{i=1}t_i^{m_i})\prod_{i=\mu}^{\nu}
(a_ie^{-\partial_i})^{m_i}f\right]\omega
=(d_j,\alpha)e_{\alpha}.f\omega.
$$
Moreover, for $\mu\le j\le \nu$, we have
$$
(d_je_{\alpha}-e_{\alpha}d_j).f\omega
=d_j.\left[(\prod^{\mu-1}_{i=1}t_i^{m_i})\prod_{i=\mu}^{\nu}
(a_ie^{-\partial_i})^{m_i}f\right]\omega-e_{\alpha}.(t_jf\omega)
$$
$$
=t_j\left[(\prod^{\mu-1}_{i=1}t_i^{m_i})\prod_{i=\mu}^{\nu}
(a_ie^{-\partial_i})^{m_i}f\right]\omega
$$
$$
-\left[(\prod^{\mu-1}_{i=1}t_i^{m_i})\prod_{i=\mu}^{\nu}
(a_ie^{-\partial_i})^{m_i}(t_jf)\right]\omega
$$
$$
=m_j\left[(\prod^{\mu-1}_{i=1}t_i^{m_i})\prod_{i=\mu}^{\nu}
(a_ie^{-\partial_i})^{m_i}f\right]\omega
=(d_j,\alpha)e_{\alpha}.f\omega
$$
as required,
where we have used the fact
$$
(a_je^{\partial_j})^{m_j}(t_jf)=(t_j-m_j)(a_je^{-\partial_j})^{m_j}f
$$
in the second to last identity. Other relations can be checked easily 
and so are omitted here. 
 
Next we prove that the module is simple.
 Let $M$ be a nonzero submodule of $M_{\omega}$. One can take a nonzero
element $f\omega\in M$ so that $f\in \C[t_1,\cdots,t_{\nu}].$
{} From the action
$d_j.f\omega=t_j(\partial_jf)\omega+f_jf\omega$, we see that 
$(\partial_jf
)\omega\in M$ for all $1\le j\le \mu-1$. This implies that one can 
choose
the nonzero polynomial $f\in \C[t_{\mu},\cdots, t_{\nu}]$, such that $f$ 
has minimum degree. We apply $e_{\alpha}$, for $\alpha\in L_C$, on
$f\omega$ to get
$$
f(t_{\mu}-m_{\mu},\cdots, t_{\nu}-m_{\nu})\omega\in M
$$
for any $m_{\mu},\cdots, m_{\nu}\in \Z$. Therefore, if $f$ depends on
$t_j$ for some $\mu\le j\le\nu$, then
$$
(f(t_{\mu},\cdots,t_j+1,\cdots,t_{\nu})-f(t_{\mu},\cdots, 
t_{\nu}))\omega\in
M
$$
and $f(t_{\mu},\cdots,t_j+1,\cdots,t_{\nu})-f(t_{\mu},\cdots, t_{\nu})$ is
nonzero with lower degree, which is a contradiction. Thus $f$ is a
nonzero constant, and so $\omega\in M$. This finishes the proof.
\qed
 
\bt{t3.4} Let $\omega_1=\omega_1(f_1,\cdots,
f_{\mu-1}|
a_{\mu},\cdots, a_{\nu})$ and $\omega_2=\omega_2(g_1,\cdots,$
 $g_{\gamma-1}|
b_{\gamma},\cdots, b_{\nu})$. Then, $M_{\omega_1}\cong M_{\omega_2}$ as 
 $B$-modules  if
and
only if $\mu=\gamma$, $a_j=b_j$ for $\mu \le j \le\nu$, and 
$f_j-g_j\in\Z$
for $1\le j\le \mu-1.$
\et
 
\pf Suppose $M_{\omega_1}\cong M_{\omega_2}$, and
$\phi:$ $M_{\omega_1}\to M_{\omega_2}$ is an isomorphism such that $
\phi(h\omega_1)=\omega_2$ for some $h\in
\C[t_1^{\pm 1},\cdots,t_{\mu-1}^{\pm 1},t_{\mu}\cdots,t_{\nu}]$.
  We first prove
$\mu=\gamma$. Otherwise one may assume $\mu>\gamma$. We have
$$
\phi(t_{\gamma}h\omega_1)=\phi(e_{c_{\gamma}}.h\omega_1)
=e_{c_{\gamma}}.\phi(h\omega_1)=e_{c_{\gamma}}.\omega_2
=b_{\gamma}\omega_2=\phi(b_{\gamma}h\omega_1)
$$
which implies that $t_{\gamma}h\omega_1=b_{\gamma}h\omega_1$, and hence
$h=0$ a contradiction. Next we prove $a_i=b_i$, for $\mu\le i\le \nu$.
Otherwise, we may assume $a_j\not= b_j$ for some $j$. Then
$$
\phi(a_jh(t_1,\cdots, t_j-1,\cdots, 
t_{\nu})\omega_1)=\phi(e_{c_j}.h\omega
_1)=e_{c_j}.\phi(h\omega
_1)
$$
$$
=e_{c_j}.\omega_2=b_j\omega_2=\phi(b_jh\omega_1)
$$
which gives
$$
a_jh(t_1,\cdots, t_j-1,\cdots, t_{\nu})\omega_1=b_jh\omega_1.
$$
This implies that $h$ is independent of $t_j$. Moreover $a_j\not= b_j$
also forces $h=0$, contradiction. Thus $a_i=b_i$ for all $i$. Finally 
we
prove
$f_j-g_j\in\Z$ for $1\le j\le \mu-1$. Set $\lambda_i=f_i-g_i$. We have
$$
\phi(t_i(\partial_ih)\omega_1+f_ih\omega_1)=\phi(d_i.h\omega_1)=
d_i.\phi(h\omega_1)
$$
$$
=d_i.\omega_2=g_i\omega_2=g_i\phi(h\omega_1).
$$
Moreover it is easy to see that 
$g_i\phi(h\omega_1)=\phi(g_ih\omega_1)$,
as $g_i\in\C[t_1^{\pm 1},\cdots, t_{\mu-1}^{\pm 1}]$. Therefore we get
$$
t_i\partial_ih=\lambda_ih
$$
where $\lambda_i=f_i-g_i$, for $1\le i\le \mu-1$.

We claim, for $1\le i\le \mu-1$, that $\lambda_i$ is independent of
$t_i$. To prove this statement, we
first introduce some notation. For $p=\sum^s_{i=r}a_it^i$ with $a_r,
a_s\not= 0$, we define deg$^+_tp=s$, deg$^-_tp=r$. If $p=0$, we define
$\deg^+_tp=0=\deg^-_tp$.

If the statement is false, then ${\hbox {deg}}^+_{t_i}\lambda_i=m^+>0$ 
or
${\hbox {deg}}^-_{t_i}\lambda_i=m^-<0$. Set ${\hbox
{deg}}^{\pm}_{t_i}h=n^{\pm}$. If $m^+> 0$ then
$$
n^+\ge {\hbox {deg}}^+_{t_i}(t_i\partial_ih)={\hbox
{deg}}^+_{t_i}\lambda_ih=m^++n^+> n^+
$$
which is a contradiction. Similarly we will get contradiction for the 
case
$m_-<0$. Thus $\lambda_i$ is independent of
$t_i$. Furthermore we claim that $\lambda_i$ is independent of
$t_j$ for any $j=1,\cdots, \mu-1$. Indeed, let $1\le j\not= i\le\mu-1$
be fixed. We know that $\lambda_i$ and $\lambda_j$ are independent of
$t_i$ and $t_j$ respectively, and also $\lambda_i,\lambda_j\in
\C[t_1^{\pm 1},\cdots, t_{\mu-1}^{\pm 1}]$, such that
$$
\partial_ih=t_i^{-1}\lambda_ih,
$$
$$
\partial_jh=t_j^{-1}\lambda_jh,
$$
where $h\in \C[t_1^{\pm 1},\cdots, t_{\mu-1}^{\pm 1},t_{\mu},\cdots,
t_{\nu}]$.

Notice that
$$
{\frac {\partial^2h}{\partial t_j\partial t_i}}={\frac 
{\partial}{\partial
t_j}}(t_i^{-1}\lambda_ih)=t_i^{-1}({\frac {\partial\lambda_i}{\partial
t_j}}h+\lambda_i{\frac {\partial h}{\partial t_j}})
=t_i^{-1}({\frac {\partial\lambda_i}{\partial
t_j}}h+\lambda_it_j^{-1}\lambda_jh),
$$
and
$$
{\frac {\partial^2h}{\partial t_i\partial t_j}}={\frac 
{\partial}{\partial
t_i}}(t_j^{-1}\lambda_jh)=t_j^{-1}({\frac {\partial\lambda_j}{\partial
t_i}}h+\lambda_j{\frac {\partial h}{\partial t_i}})
=t_j^{-1}({\frac {\partial\lambda_j}{\partial
t_i}}h+\lambda_jt_i^{-1}\lambda_ih),
$$
this gives us
$$
(t_i^{-1}{\frac {\partial\lambda_i}{\partial
t_j}}-
t_j^{-1}{\frac {\partial\lambda_j}{\partial
t_i}})h=0
$$
or
$$
t_j{\frac {\partial\lambda_i}{\partial
t_j}}-
t_i{\frac {\partial\lambda_j}{\partial
t_i}}=0
$$
as $h\not= 0$.
Set $\Phi=t_i{\frac {\partial\lambda_j}{\partial
t_i}}\in \C[t_1^{\pm 1},\cdots,t_{\mu-1}^{\pm 1}]$, which is 
independent
of
 $t_j$. Then
 $
t_j{\frac {\partial\lambda_i}{\partial
t_j}}=\Phi
$
 gives us $\lambda_i=\Phi\ln t_j+C\in
\C[t_1^{\pm 1},\cdots,t_{\mu-1}^{\pm 1}]$.
Thus $\Phi=0$. That is $t_j{\frac {\partial\lambda_i}{\partial
t_j}}=0$, or $\lambda_i$ is independent of $t_j$. This proves that
$\lambda_i\in\C$ for $i=1,\cdots,\mu-1$. Moreover the 
equation
$$
{\frac {\partial h}{\partial t_i}}=t^{-1}_i\lambda_ih
$$
gives $h=Ct^{\lambda_i}_i$, but $h\in
\C[t_1^{\pm 1},\cdots,t_{\mu-1}^{\pm 1},t_{\mu},\cdots, t_{\nu}]$ 
forces
$\lambda_i\in\Z$ as $h\not= 0$.

 Now we suppose the conditions  hold, and prove the two modules
$M_{\omega_1}$ and $M_{\omega_2}$ are isomorphic, where
$$
\omega_1=\omega_1(f_1,\cdots, f_{\mu-1}|a_{\mu},\cdots, a_{\nu})
$$
$$
\omega_2=\omega_2(g_1,\cdots, g_{\mu-1}|a_{\mu},\cdots, a_{\nu})
$$
and $g_i=f_i+N_i$ for some $N_i\in\Z$. We define $\phi:$ $M_{\omega_1}
\to M_{\omega_2}$ by
$$
\phi(f\omega_1)=t_1^{-N_1}\cdots t_{\mu-1}^{-N_{\mu-1}}f\omega_2
$$
for $f\in \C[t_1^{\pm 1},\cdots,t_{\mu-1}^{\pm 1},t_{\mu},\cdots,.
t_{\nu}].$
 It is obvious that $\phi$ defines a vector space isomorphism. We need
to check that it also defines an algebra homomorphism.

Let $\alpha=\sum^{\nu}_{i=1}m_ic_i$. It is easy to see that
$$
\phi(e_{\alpha}.f\omega_1)=e_{\alpha}\phi(f\omega_1),
$$
and $\phi(d_l.f\omega_1)=d_l.\phi(f\omega_1)$ for $\mu\le l\le\nu$.
Therefore we only need to check the identity
$\phi(d_l.f\omega_1)=d_l.\phi(f\omega_1)$ for $1\le l\le\mu-1$. In fact
$$
d_l.\phi(f\omega_1)=d_l.(\prod^{\mu-1}_{i=1}t_i^{-N_i})f\omega_2
$$
$$
=\left[ t_l\partial_l\left( (\prod^{\mu-1}_{i=1}t_i^{-N_i})f\right)+
g_l(\prod^{\mu-1}_{i=1}t_i^{-N_i})f\right]\omega_2
$$
$$
=\left[
-N_l(\prod^{\mu-1}_{i=1}t_i^{-N_i})f+t_l(\prod^{\mu-1}_{i=1}t_i^{-N_i})
\partial_lf+g_l(\prod^{\mu-1}_{i=1}t_i^{-N_i})f\right]\omega_2
$$
$$
=\left[
(\prod^{\mu-1}_{i=1}t_i^{-N_i})t_l\partial_lf+(\prod^{\mu-1}_{i=1}t_i^{-N_i})
f_lf\right]\omega_2
$$
$$
=\phi((t_l\partial_lf+f_lf)\omega_1)=\phi(d_lf\omega_1),
$$
as required. This finishes the proof of the Theorem.
\qed

Now we  define the associative algebra $A$ to be the  quotient
of $B$ modulo relations $d_id_j=d_jd_i$ for all $i,j.$ 
\begin{theorem}\label{t3.5}
For $\omega=\omega(f_1,\cdots, f_{\mu-1}|a_{\mu},\cdots,
a_{\nu}),$ $M_{\omega}$ is an A-module if and only if $f_j=t_j\partial_jP+P_j(t_j) (
j=1,\cdots,\mu-1)$ for some $P\in \C[t_1^{\pm 1},\cdots, t_{\mu-1}^{\pm 1}]$
and $P_j(t_j)\in \C[t_j^{\pm 1}]$.
\end{theorem}
 
\pf  The proof of the "if" part is straightforward. To prove the "only if" 
part, we
suppose $M_{\omega}$ is an A-module. From the relation $d_id_j=d_jd_i$ for $
1\le i,j\le \mu-1$, one can easily obtain $D_if_j=D_jf_i$ for $1\le i,j\le
\mu-1$, where $D_j=t_j\partial_j$ is the degree derivation of $\C[t_1^{\pm 1},
\cdots, t_{\mu-1}^{\pm 1}]$.
 
  For this fixed $\omega=\omega(f_1,\cdots,f_{\mu-1}|a_{\mu},\cdots,a_{\nu})$,
we can write $f_1,\cdots, f_{\mu-1}$ in the following form
$$
f_i=\sum_{k_1,\cdots,k_{\mu-1}\in Z}
a_{k_1,\cdots,k_{\mu-1}}^{(i)}t_1^{k_1}\cdots
t_{\mu-1}^{k_{\mu-1}}+P_i(t_i)
$$
for $i=1,2,\cdots,\mu-1$, where $
a_{k_1,\cdots,k_{\mu-1}}^{(i)}\in \C$, and the homogeneous term $
t_1^{k_1}\cdots t_{\mu-1}^{k_{\mu-1}}\not\in \C[t_i^{\pm 1}]$ if $
a_{k_1,\cdots,k_{\mu-1}}^{(i)}\not= 0$.
 
Suppose $a_{k_1,\cdots, k_{\mu-1}}^{(i)}\not= 0$ for a fixed $i$, then, from $
D_j f_i=D_if_j$, we have $k_j
a_{k_1,\cdots, k_{\mu-1}}^{(i)}=k_ia_{k_1,\cdots, k_{\mu-1}}^{(j)}$, which implies
that $k_i\not= 0$(as $t_1^{k_1}\cdots t_{\mu-1}^{k_{\mu-1}}\not\in \C[t_i^{
\pm 1}]$, thus one can find $j(\not= i)$ so that $k_j\not= 0$). Let $
c_{k_1,\cdots, k_{\mu-1}}={\frac 1 {k_i}}a_{k_1,\cdots, k_{\mu-1}}^{(i)}\in \C$.
It follows from $k_ja_{k_1,\cdots, k_{\mu-1}}^{(i)}=k_i
a_{k_1,\cdots, k_{\mu-1}}^{(j)}$ that
$$
a_{k_1,\cdots, k_{\mu-1}}^{(j)}=c_{k_1,\cdots, k_{\mu-1}}k_j
$$
for all $1\le j\le \mu-1$. Therefore
$$
f_j=D_j(\sum_{k_1,\cdots, k_{\mu-1}\in \Z}c_{k_1,\cdots, k_{\mu-1}}
t_1^{k_1}\cdots t_{\mu-1}^{k_{\mu-1}})+P_j(t_j)
$$
for $1\le j\le \mu-1$, as required. This also completes the proof. \qed

\begin{corollary}\label{c3.6} For $\omega\!=\!\omega(\!f_1,\cdots, f_{\mu-1}|a_{\mu},\cdots, a_{\nu}\!)$, $M_{\omega}
$ is a simple A-module if and only if $f_j=
D_jP+P_j(t_j)$ $(j=1,\cdots, \mu-1)$ for some Laurent polynomials $P\in
\C[t_1^{\pm 1},\cdots, t_{\mu-1}^{\pm 1}]$ and $P_j(t_j)\in \C[t_j^{\pm 1}]$.
Moreover, if $M_{\omega_1}$ and $M_{\omega_2}$ are A-modules with $\omega_1=
\omega_1(f_1,\cdots, f_{\mu-1}|a_{\mu},\cdots, a_{\nu})$ and $\omega_2=
\omega_2(g_1,\cdots, g_{\gamma-1}|b_{\gamma},\cdots, b_{\nu})$, then $
M_{\omega_1}\cong M_{\omega_2}$ as A-modules if and only if $\mu=\gamma$, $
a_j=b_j$ $(\mu\le j\le \nu)$, and $f_j-g_j\in \Z$ $(1\le j\le \mu)$.
\end{corollary}
 
For $\alpha=\sum^{\nu}_{i=1}m_ic_i\in L_C$, we
set $t^{\alpha}=t_1^{m_1}\cdots t_{\nu}^{m_{\nu}}$ $\in$ $\C[t_1^{\pm
1},\cdots, t_{\nu}^{\pm 1}].$ Let $\mu=\nu$ in the previous theorem, 
and
take $f_i=\lambda_i\in\C$ for $1\le i\le \nu$. Then 
$M_{\omega}=\oplus_{
\alpha\in L_C}\C t^{\alpha}\omega$, where
$\omega=\omega(\lambda_1,\cdots,
\lambda_{\nu}|)$. It is clear that
$$
d_i.(t^{\alpha}\omega)=(d_i,\alpha+\lambda)t^{\alpha}\omega
$$
for $1\le i\le\nu$, where $\lambda=\sum^{\nu}_{i=1}\lambda_ic_i\in 
\h_C.$
 Thus we have
\begin{corollary}\label{c3.7} If $\omega=\omega(\lambda_1,\cdots,
\lambda_{\nu}|)$ with $\lambda=\sum^{\nu}_{i=1}\lambda_ic_i\in H_C$ 
then $M_{\omega}$ is isomorphic to $ \C[L_C+\lambda]$ as $A$-modules.
Moreover, $\C[L_C+\lambda]\cong\C[L_C+\lambda']$ if and only if $
\lambda-\lambda'\in L_C.$
\end{corollary}

The connection between $A$-modules and $V$-modules is studied in the 
next two
sections.

\section{Construction\,of\,$V$-Modules\,from\,$A$-Modules}
Let $W$ be an $A$-module. For $\lambda\in{\frac 1 k}L_D$, set
$$
V_{\lambda,W}=M(1)\otimes W.
$$
Our objective is to make $V_{\lambda,W}$ a $V$-module.

Motivated by the representation theory for the lattice vertex algebras,
we define an action of $\hat{\frak h}$ on $V_{\lambda,W}$ as follows:
\begin{eqnarray*}
& &\alpha(n)\mapsto \alpha(n)\otimes 1\\
& &\beta(0)\mapsto (\beta,\lambda)\\
& &\gamma(0)\mapsto 1\otimes \gamma\\
& &\ \ \ \ \  \ c\mapsto 1
\end{eqnarray*}
for $n\in\Z,$ $\alpha\in {\frak h},$
$\beta\in {\frak h}_C$, $\gamma\in {\frak h}_D.$
Then $V_{\lambda,W}$ is an $\hhh$-module. We also define
operators $e^{\alpha},$  on $V_{\lambda,W}$ and $z^{\alpha}$ on $V_{\lambda,W}[[z,z^{-1}]]$
for $\alpha\in L_C$ as follows:
$$
e^{\alpha}\mapsto 1\otimes e_{\alpha},\;\;\;\;\;\;
z^{\alpha}\mapsto z^{(\alpha,\lambda)}.
$$

\bl{l3.1} For $\alpha, \beta\in L_C$, $\gamma\in {\frak h},$ we have
$$
e^{\alpha}z^{\beta}=z^{\beta}e^{\alpha},\;\;\;\;
\gamma(0)e^{\alpha}=e^{\alpha}\gamma(0)+(\gamma,\alpha)e^{\alpha}.$$
\el

For $\alpha\in\fh$ we set
$$\alpha(z)=\sum_{n\in\zzz}\alpha(n)z^{-n-1}.$$
We also define
$$ Y_{\lambda,W}(e^{\beta},z)=E^-(-\beta,z)E^+(-\b,z)e^{\b}z^{\b}$$
for $\beta\in L_C$ where
$$E^{\pm}(\beta,z)=\exp(\sum_{n\in\pm \N}\frac{\beta(n)z^{-n}}{n}).$$
Now we define a linear map
\begin{eqnarray*}
&& Y_{\lambda, W}:\;\;\;V\to ({\hbox {End}}V_{\lambda,W})[[z,z^{-1}]]\\
&&\ \ \ \ \ \  \;\;\;\;\;\;\;\;v\mapsto 
Y_{\lambda,W}(v,z)=\sum_{n\in\zzz}v_nz^{-n-1}
\end{eqnarray*}
by
$$
Y_{\lambda,W}(\alpha_1(-n_1)\cdots \alpha_s(-n_s)e^{\alpha},z)
$$
$$
=:(\partial_{n_1-1}\alpha_1(z))\cdots
\partial_{n_s-1}\alpha_s(z))Y_{\lambda,W}(e^{\alpha},z):
$$
where the normal ordering is defined in [FLM],  $\partial_n={\frac 1
{n!}}({\frac {\partial}{\partial z}})^n.$

\br{flm} {\rm The expression of the operator $Y_{\lambda,W}(v,z)$ here
is the same as  that in (\ref{5.8}) except that we
deal with abstract space $W.$}
\er

It is easy to see that $V_{\lambda, W}$ is a module for the Heisenberg
vertex algebra $M(1)$. The rest of this section is devoted to the proof 
that
$V_{\lambda, W}$ is a $V$-module. As we have already pointed out that  
we cannot prove the Jacobi identity 
for the operators $Y_{\lambda,W}(u,z)$ by using the
coherent state argument given in Chapter 8 of [FLM] as our lattice
$L_C$ does not {\em span} $\fh.$ Instead we will apply the theory
of local system developed in [L2] to our situation. 
 
Recall that $\omega=\frac{1}{2}\sum_{i=1}^{2\nu}\beta_i(-1)^2.$ It 
follows from
the definition of $Y_{\lambda,W}(v.z)$ that
\begin{equation}\label{add1}
Y_{\lambda,W}(L(-1)v,z)=\frac{d}{dz}Y_{\lambda,W}(L(-1)v,z)
\end{equation}
for all $v\in V.$

We need several lemmas.

Recall (\ref{action}). The following lemma is straightforward.

\bl{addl1} Let $\a_i\in\fh, 0<n_i\in \Z$ for $i=1,...,s$ and
$\a\in L_C.$ Then
$$Y_{\lambda,W}(\alpha_1(-n_1)\cdots \alpha_s(-n_s)e^{\alpha},x)
=\a_1(x)_{-n_1}\cdots\a_s(x)_{-n_s}Y(e^{\a},x).$$
\el

Set $S_{\lambda,W}=\{Y_{\lambda,W}(v.x)|v\in V\}.$ We have the 
following:

\bl{l3.2} $S_{\lambda,W}$ is a local space on $V_{\lambda,W}.$
\el

\pf Let $h,h'\in{\frak h}$ and $\alpha,\beta\in L_C.$
Since $V_{\lambda,W}$ is a module for the affine algebra
$\hhh$ we immediately have
$$[h(x_1),h'(x_2)](x_1-x_2)^2=0.$$
One can easily verify that
$$[h(x_1),Y_{\lambda,W}(e^{\beta},x_2)](x_1-x_2)=0$$
$$[Y_{\lambda,W}(e^{\alpha},x_1),Y_{\lambda,W}(e^{\beta},x_2)]=0$$
(cf. Chapter 7 of [FLM]). Thus
$X=\<\alpha(x), Y_{\lambda,W}(e^{\beta},x)|\alpha\in {\frak h},\beta\in 
L_C\>$
is a local space on  $V_{\lambda,W}.$
By Lemma \ref{addl1} $S_{\lambda,W}$
can be generated from $X$ by the operations $a(x)_nb(x)$ and by tacking derivations.
Thus $S_{\lambda,W}$ is a local space on $V_{\lambda,W}$ (see [L2]).
\qed

\br{addr1} {\rm Lemmas \ref{addl1} and \ref{l3.2} also hold
for vertex operators $Y(u,z)$ acting on $V$ (see [FLM], [DL]).\rm}
\er

Now we set $M=V\oplus V_{\lambda,W}$ and define
$u(z)=Y_M(u,z)=Y(u,z)+Y_{\lambda,W}(u,z)$ for $u\in V.$ We also set
$S=\<Y_M(u,x)|u\in V\>,$
a local space on $M.$
Let $U$ be a vertex algebra generated by the local
space $S.$  By Theorem \ref{tli}, $M$ is an
$U$-module with
$Y_M(a(x),z)=a(z)$ for $a(x)\in U$ and $Y_M(u(x),z)=Y_M(u,z)$ if $u\in 
V$.
 
By Lemma 2.2.5 of [L2], $U$ is a module for the affine algebra
$\hhh$ with $h(n)$ acting as $h(x)_n$ for $h\in \fh,$ $n\in\Z$ and $c$
acting as 1. Thus $U$ contains a vertex subalgebra $M(1)$ generated by
$h(x)$ for $h\in \fh.$ We should mention that the component
operators $\omega(x)_n$ for $n\in\Z$ satisfy the Virasoro
algebra relation
$$[\omega(x)_{m+1},\omega(x)_{n+1}]=(m-n)\omega(x)_{m+n+1}+
\frac{m^3-m}{6}\delta_{m+n,0}\nu$$
which follows directly from the affine algebra relation in $\hhh.$

Recall the component operators of $Y(u(x),z)$ are given by
$Y(u(x),z)=\sum_{n\in\Z}u(x)_nz^{-n-1}$
where $Y$ is the linear map from $U$ to $(\End U)[[z,z^{-1}]]$ which is guaranteed 
 since $(U,Y,I(x),\frac{d}{dx})$ is a vertex algebra.

\bl{addl3.2} Let $h\in\fh, m,n\in\Z, \a\in L_C.$ Then
\begin{equation}\label{adde1}
[h(x)_n,e^{\alpha}(x)_m]=(h,\a)e^{\a}(x)_{m+n}.
\end{equation}
\el
\pf Since $U$ is a vertex algebra we have the commutator formula
which is a consequence of the Jacobi identity:
$$[h(x)_n,e^{\alpha}(x)_m]=\sum_{i\geq 0}{n\choose 
i}(h(x)_ie^{\a}(x))_{m+n-i}.$$
Note that
\begin{eqnarray*}
& &h(x)_0e^{\alpha}(x)=\res_z\{h(z)e^{\alpha}(x)-e^{\alpha}(x)h(z)\}\\
& &\ \ \ \ \ \ \ \ =[h(0),Y(e^{\alpha},x)]=(h,\alpha)Y(e^{\alpha},x).
\end{eqnarray*}
If $i\geq 1$ we know from the proof of Lemma \ref{l3.2} that
$$h(x)_ie^{\alpha}(x)=\res_z\{(z-x)^ih(z)e^{\alpha}(x)-(z-x)^ie^{\alpha}(x)h(z)\}=0.$$
The result follows immediately. \qed

\bl{addl3.3} Let $\a,\b\in L_C,$ and $m\in \Z$  Then
$e^{\a}(x)_me^{\b}(x)=0$ if $m\geq 0$ and
$e^{\a}(x)_{m}e^{\b}(x)=Y(u^m\otimes e^{\a+\b},x)$ if
$m<0$ where $u^m\in M(1)$ is determined by 
$\frac{1}{(-m-1)!}L(-1)^{-m-1}e^{\a}=e^{\alpha}_m=u^m\otimes e^{\a}.$
\el

\pf From the proof of Lemma \ref{l3.2} we know that if $m\geq 0$ then
$$[e^{\a}(z),e^{\b}(x)](z-x)^m=0$$
 of nonnegative $m.$ It is immediate
{}from the definition of $e^{\a}(x)_me^{\beta}(x)$ that
$e^{\a}(x)_me^{\b}(x)=0$ in this case.

Now we deal with negative $m.$ If $m=-1$ then
$$e^{\a}(x)_{-1}e^{\b}(x)=Y(e^{\a},x)^{-}Y(e^{\b},x)+
Y(e^{\b},x)Y(e^{\\a},x)^{+}$$
where
$$Y(e^{\\a},x)^{+}=\sum_{s\geq 0}e^{\alpha}_sx^{-s-1},
Y(e^{\\a},x)^{-}=\sum_{s<0}e^{\alpha}_sx^{-s-1}.$$
Since
$$[e^{\alpha}_s,e^{\beta}_t]=0$$
for any $s,t\in\Z$ and $(\alpha,\beta)=0,$  we have
$$e^{\a}(x)_{-1}e^{\b}(x)=Y(e^{\alpha},x)Y(e^{\b},x)=Y(e^{\alpha+\beta},x).$$
Now let $m=-1-n$ for some nonnegative $n.$ A straightforward computation
using (\ref{add1}) gives

\begin{eqnarray*}
& 
&e^{\a}(x)_{m}e^{\b}(x)=\frac{1}{n!}\left(\frac{d}{dx})^nY(e^{\a},x)\right)
Y(e^{\b},x)\\
& &\ \ \ \ =\frac{1}{n!}Y(L(-1)^ne^{\a},x)Y(e^{\b},x)\\
& &\ \ \ \ =Y(u^m\otimes e^{\a+\b},x),
\end{eqnarray*}
as required. \qed

\bp{padd1} We have $S=U.$ That is, $S$ is a vertex algebra.
\ep

\pf Since $U$ is generated by
$$X=\<h(x), e^{\alpha}(x)|h\in \fh, \a\in L_C\>$$
{}from the proof of Lemma \ref{l3.2} and Remark \ref{addr1}, it is 
enough
to show that $S$ in invariant under the operator $u(x)_n$ for
any $u(x)\in X$ and $n\in \Z.$

Note from Lemma \ref{addl1} that $S$ is spanned by
$$h_1(x)_{-n_1}\cdots h_s(x)_{-n_s}e^{\b}(x)$$
for $h_i\in \fh,$ $n_i>0$ and $\beta\in L_C.$ Also note that
$h(x)_ne^{\b}(x)=(h,\b)\delta_{n,0}e^{\b}(x)$ for nonnegative
$n$ (see the proof of Lemma \ref{addl3.2}). Using the commutator
relation $[h(x)_m,h'(x)_n]=m\delta_{m,-n}(h,h')$ we see that
$S$ is invariant under $h(x)_n$ for all $h\in\fh, n\in \Z.$

Let $\a\in L_C.$ By Lemma \ref{addl3.3}, $e^{\a}(x)_ne^{\b}(x)\in S.$
An induction by using Lemma \ref{addl3.2} on $s$ then shows
that  $e^{\a}(x)_nh_1(x)_{-n_1}\cdots h_s(x)_{-n_s}e^{\b}(x)\in S.$
\qed

We are now in a position to prove the main result in this section.

\bt{addt1} Let $\lambda, W$ be as before. Then $(V_{\lambda,W}, 
Y_{\lambda,W})$
is a $V$-module. Moreover, $V_{\lambda,W}$ is irreducible if and only
if $W$ is simple.
\et

\pf By Theorem \ref{tli} and Proposition \ref{padd1}, $V_{\lambda,W}$
is a module for the vertex algebra $S$ under the action
$Y_{V_{\lambda,W}}(u(x),z)=u(x)$ for any $u(x)\in S.$ Note that
$u(x)$ on $V_{\lambda,W}$  is exactly $Y_{\lambda,W}(u,x)$ for $u\in 
V.$
So it is enough to prove that the map from $V$ to $S$ 
 defined by sending $u$ to $Y(u,x)=u(x)$ is a vertex algebra isomorphism.

Let $\bar V=\{Y(u,x)|u\in V\}$ be the set of vertex operators on $V.$ Then
$\bar V$ is a vertex algebra under
$$Y(u,x)_nY(v,x)=\Res_\{(z-x)^nY(u,z)Y(v,x)-(-x+z)^nY(v,x)Y(u,z)\}$$
for $u,v\in V$ and both maps from $V$ and $S$ to $\bar V$ given
by sending $v$ and $v(x)$ to $Y(v,z)$ for $v\in V$ are
surjective vertex algebra homomorphisms (see [L2]). Note that
$Y(v,z)=0$ if and only if $v=0.$ This shows that both homomorphisms
are isomorphisms. As a result, the map from $V$ to $S$ given by sending
$v$ to $v(x)$ is a vertex algebra isomorphism.

Now we assume that $V_{\l,W}$ is an irreducible $V$-module. Let $M$ a nonzero
$A$-submodule. It is clear from the definition of vertex operators
$Y_{\l,W}(v,z)$ that $M(1)\otimes M$ is a $V$-submodule
of $V_{\l,W}.$ Thus $W$ is simple. It will be proved in the
next section that if $W$ is a simple $A$-module then $V_{\l,W}$
is an irreducible $V$-module. \qed
  
\section{Construction\,of\,$A$-Modules\,from\,$V$-Modules}

Let $M=(M,Y_M)$ be a $V$-module. Set 
$Y_M(\alpha(-1),z)=\sum_{n\in\Z}\alpha(-n)z^{-n-1}$ for $\alpha\in \fh.$ Then the operators $\alpha(m), 
\beta(n)$
satisfy the Heisenberg algebra relation
\begin{equation}\label{4.1}
[\alpha(m),\beta(n)]=m(\alpha,\beta)\delta_{m+n,0}
\end{equation}
(see [D] for the details).
Thus $M$ is a module for the affine algebra $\hat\fh.$ We now define
the vacuum space $\Omega_M$ following [LW] (also see [D]):
$$\Omega_M=\{w\in M|\a(n)w=0, \a\in \fh, n>0\}.$$
In this section we prove that ``weight space'' $\Omega_M'$ of
$\Omega_M$ is an $A$-module.  We should point out that it is possible 
that
$\Omega_M=0$ without additional assumption.

We set $\bar M=M(1)\otimes \Omega_M$  then $\bar M$ is a subspace
of $M$ and we will also  show that $\bar M$ is a $V$-submodule of $M.$

For $\alpha\in L_C$, set
$$
Z(\alpha,z)=\exp\left( \sum_{n\in\zzz_-}{\frac
{\alpha(n)}{n}}z^{-n}\right)Y_M(e^{\alpha},z)\exp\left(
\sum_{n\in\zzz_+}{\frac
{\alpha(n)}{n}}z^{-n}\right)
$$
$$
=\sum_{n\in\zzz}Z(\alpha, n)z^{-n-1}.
$$
\bl{l4.1} For $\alpha\in L_C,$ $\beta\in H$ we have

(1) $[\beta(0),Z(\alpha, z)]=(\beta,\alpha)Z(\alpha,z),$

(2) $[\beta(n),Z(\alpha, z)]=0,$ if $n\not= 0$,

(3) ${\frac {d}{dz}}Z(\alpha,z)=Z(\alpha, z)\alpha(0)z^{-1},$

(4) $Z(\alpha,n)\Omega_M\subset \Omega_M$ for $n\in\zzz.$

(5) $[\alpha(0), Y_M(v,z)]=0,$ for $v\in V$.
\el

\pf The proofs follow from a similar argument as in
[D]. We refer the reader to [D] for details. \qed

\bl{l4.2} If $\Omega_M\ne 0$ there exists a nonzero
$w\in \Omega_M$ and $\lambda\in \frac{1}{k}L_D$ such that
$\alpha(0)w=(\lambda,\alpha)w$ for all $\alpha\in L_C.$
\el

\pf One of the important consequence of the Jacobi identity (\ref{2.2})
is  associativity (see [DL], [FLM], [L2]): Let $u,v\in V$ and
$w\in M,$ then there exists $n\geq 0$ such that
$$(z_2+z_0)^nY_M(u,z_0+z_2)Y_M(v,z_2)w=Y_M(Y(u,z_0)v,z_2)w.$$
Using this associativity one can show (see the proof of
Corollary 4.2 of [DM] or the proof of Proposition 4.1 of [L1]),
 that if $Y_M(u,z)w=0$ for some nonzero $u\in V$ and
$w\in W$ then $Y(v,z)w=0$ for all $v\in V.$ (One needs to use
the fact that $V$ is simple.) Since $Y_M({\bf 1},z)$ is
the identity operator, we conclude that for any nonzero $u\in V$
and nonzero $w\in M,$ $Y_M(u,z)w\ne 0.$  In particular,
$Y_M(e^{\alpha},z)w\ne 0$ for any $w\in M.$

{}From the definition of operator $Z(\alpha,z)$ we know that
$$Y_M(e^{\alpha},z)=\exp\left( \sum_{n\in\zzz_-}{-\frac
{\alpha(n)}{n}}z^{-n}\right)\exp\left(\sum_{n\in\zzz_+}{-\frac
{\alpha(n)}{n}}z^{-n}\right)\otimes Z(\alpha,z)$$
on $M(1)\otimes \Omega_M$ by Lemma \ref{l4.1} (4).
Thus $Z(\alpha,z)w\ne 0$ for any nonzero $w\in \Omega_M.$
By Lemma \ref{l4.1} (3) we have the following
$$\alpha(0)Z(\alpha,n)=(-n-1)Z(\alpha,n)$$
for any $n\in\Z.$  The exact same proof of Lemma 3.4 in [D] then
works here. \qed

Motivated by Lemma \ref{l4.2} we define a weight subspace
$$\Omega_M(\lambda)=\{w\in\Omega_M|\alpha(0)w=(\alpha,\l)w\ {\rm for\ 
all \ }
\alpha\in L_C\}$$
for any $\l\in L_D.$  Set $\Omega_M'=\sum_{\l\in L_D}\Omega_M(\l).$
Now we define an operator $z^{\alpha}$ on $\Omega_M'$ for $\alpha\in L_C$ 
 by saying it acts  on $\Omega_M(\l)$ as $z^{(\alpha,\l)}.$ Set
$$T_{\alpha}=Z(\alpha,z)z^{-\alpha}.$$
Then by Lemma \ref{l4.1} (3),
$$\frac{d}{dz}T_{\a}=0.$$
That is, $T_{\a}$ is an operator on $\Omega_M'$ and is independent
of the formal variable $z.$

\bl{l4.4} For $\alpha,\beta\in L_C, d\in \fh,$ we have
\begin{equation}\label{e4.1}
[d(0),T_{\alpha}]=(d,\alpha)T_{\alpha},\;\;\;\; 
T_{\alpha}z^{\beta}=z^{\beta}T_{\alpha},
\end{equation}
and
\begin{equation}\label{e4.2}
T_{\alpha}T_{\beta}=T_{\alpha+\beta}.
\end{equation}
\el

\pf (\ref{e4.1}) follows
{}from Lemma \ref{l4.1} by noting that $d(0)z^{\a}=z^{\a}d(0).$

For (\ref{e4.2}) we note that
$$
Y_M(e^{\alpha}, z)=E^-(-\alpha, z)E^+(-\alpha, z)T_{\alpha}z^{\alpha}.
$$
{}From the proof of Lemma \ref{addl3.3},
$$Y(e^{\alpha+\beta},z)=
Y_M(e^{\alpha}, z)Y_M(e^{\beta}, z)=E^-(-\alpha-\beta, z)E^+(-\alpha
-\beta,
z)T_{\alpha}T_{\beta}z^{\alpha+\beta}.$$
Thus we obtain $T_{\alpha}T_{\beta}=T_{\alpha+\beta}.$ \qed

\bt{t4.1} Let $M$ be a $V$-module.  Then

(1) $\bar M=M(1)\otimes \Omega_M$ is a $V$-submodule of $M.$

(2) $\Omega'_M$ is an $A$-module
by sending $d$ to $d(0)$ and $e_{\alpha}$ to $T_{\alpha}$
for $d\in H_D$ and $\alpha\in L_C$.

(3) If $M$ is irreducible and $\Omega_M\ne 0$ then
$\Omega_M'=\Omega_M$ and $M$ is isomorphic to $V_{\l,\Omega_M}$
for some $\l\in L_D.$
\et
 
\pf Clearly, $\bar M$ is invariant under the vertex operators
$Y_M(d(-1),z)$ and $Y_M(e^{\alpha},z)$ for $d\in \fh$ and $\alpha\in 
L_C.$
We have already mentioned that $V$ is generated by $d(-1)$
and $e^{\alpha}.$ It is clear now that $\bar M$ is a submodule. This
proves (1).

(2) is a direct consequence of Lemma \ref{l4.4}.

We now prove (3). By Lemma \ref{l4.2}, $\Omega_M'\ne 0.$  Set
$M'=M(1)\otimes \Omega_M'.$ Clearly $M'$ is invariant under the
component operators of $Y_M(d(-1),z)$ and $Y_M(e^{\alpha},z)$
for $d\in\fh$ and $\alpha\in L_C.$ Again since $V$ is generated by
$d(-1)$ and $e^{\alpha}$ for $d\in\fh$ and $\alpha\in L_C,$
$M'$ is invariant under the component operators of $Y_M(v,z)$ for
all $v\in V.$ That is, $M'$ is a nonzero submodule of $M.$ The
irreducibility of $M$ then yields that $M=M'.$ As a result
we have $\Omega_M=\Omega_M'.$

Note that $\alpha(0)d(n)=d(n)\alpha(0)$ and
 $\alpha(0)T_{\beta}=T_{\beta}\alpha(0)$ for all $d\in \fh$ and 
$\a,\b\in L_C$ (see Lemma \ref{l4.4}). Thus $\alpha(0)$ commutes with
$Y_M(u,z)$ for all $u\in V.$ Thus if $\alpha(0)$ for $\alpha\in L_C$ 
have
a common eigenvector $w$ with eigenvalue $\l\in L_D$ then
$\alpha(0)$ acts on $M$ as $(\alpha,\l).$ By Lemma \ref{l4.2}
we see that $V$ is isomorphic to   $V_{\l,\Omega_M}$
for some $\l\in L_D.$ \qed
 
We can now finish the proof of Theorem \ref{addt1}.
That is, if $W$ is a simple $A$-module then 
$V_{\l,W}$ is an irreducible $V$-module. Note that $\Omega_{V_{\l,W}}
=W.$ From the proof of Theorem \ref{t4.1} any $V$-submodule $M$ has 
decomposition $M=M(1)\otimes W'$ and $W'$ is an $A$-submodule of $W.$
The simplicity then implies $W'=W$ and $M=V_{\l,W}$ if $M\ne 0.$ 
 
\section{Zhu algebra $A(V)$}

In this section we compute the Zhu algebra $A(V)$ [Z] and study the 
representation
theory of $V$ in terms of $A(V).$
 
Recall from [Z] that $A(V)=V/O(V)$ where $O(V)$ is spanned by
$$u\circ v={\rm Res}_z
\frac{(1+z)^{\wt u}}{z^2}Y(u,z)v= \sum_{i=0}^\infty\binom{\wt(u)}{i}
u_{i-2}v$$
for homogeneous $u,v\in V.$ The product is defined by
$$u*v={\rm Res}_z\frac{(1+z)^{\wt u}}{z}Y(u,z)v
 =\sum_{i=0}^{\infty}{\wt u\choose i}u_{i-1}v.$$

The following lemma can be found in [Z].
\bl{l6.1} Assume that $u\in V$ homogeneous, $v\in V$ and $ n \geq 0$.
Then 
$$\Res_z\frac{(1+z)^{\wt u}}{z^{2+n}}Y(u,z)v
= \sum_{i=1}^\infty \binom{\wt u}{i}u_{i-n-2}v
\in O(V).$$
\el
 
We will write $u\sim v$ if $u-v\in O(V)$ for $u,v\in V.$ Note that
$Y(h(-1),z)=\sum_{n\in\Z}h(-1)_nz^{-n-1}=\sum_{n\in\Z}h(n)z^{-n-1}$
and the weight of $h(-1)$  is 1 for $h\in {\frak h}.$ By Lemma
\ref{l6.1} we see that
$h(-n-1)w\sim -h(-n)w$ for any $w\in V.$ Thus $A(V)$ is
spanned by
$$\C[h(-1)|h\in {\frak h}]\otimes \C[L_C].$$

For any $\alpha\in L_C$ the weight of $e^{\alpha}$ is 0. Thus
\begin{eqnarray*}
& &e^{\alpha}\circ e^{\beta}={\rm 
Res}_z\frac{1}{z^2}Y(e^{\alpha},z)e^{\beta}\\
& &\ \ \ \ ={\rm Res}_z\frac{1}{z^2}E(-\alpha,z)e^{\alpha+\beta}\\
& &\ \ \ \ =\alpha(-1)e^{\alpha+\beta}
\end{eqnarray*}
for all $\alpha,\beta\in L_C.$ Since $\alpha,\beta\in L_C$ are 
arbitrary
we conclude that $\alpha(-1)e^{\beta}\in O(V).$ In particular,
$\alpha(-1)\in O(V).$ Thus $\alpha(-1)*w=\alpha(-1)w$ lies in
$O(V)$ for any $w\in V.$ As a result we see that
$A(V)$ is spanned by
$$\C[d_1(-1),..., d_{\nu}(-1)]\otimes \C[L_C].$$  
\bp{p6.1} The Zhu algebra $A(V)$ is isomorphic to $A.$ 
\ep
 
\pf We define a linear map
$I$ from $\bar A$ to $A(V)$ by sending $d_1^{n_1}\cdots 
d_{\nu}^{n_{\nu}}e^{\alpha}$ to $d_1(-1)^{n_1}\cdots 
d_{\nu}(-1)^{n_{\nu}}e^{\alpha}+O(V).$ It is
enough to prove that $I$ is an algebra isomorphism.
 
It is a straightforward verification that $d_i(-1), d_j(-1), 
e^{\alpha},
e^{\beta}$ satisfy the relations 
\begin{equation}\label{6.2}
\begin{array}{l}
 d_i(-1)*d_j(-1)=d_j(-1)*d_i(-1)\\
d_i(-1)*e^{\alpha}-e^{\alpha}*d_i(-1)+(d_i,\alpha)e^{\alpha}\\
e^{\alpha}*e^{\beta}=e^{\alpha+\beta}
\end{array}
\end{equation}
for $1\leq i,j\leq \nu$ and $\alpha,\beta\in L_C.$
So $I$ is an onto algebra homomorphism. 
 
Proving that $I$ is injective is equivalent to proving that the 
intersection
of $O(V)$ with the subspace $\C[d_1(-1),..., d_{\nu}(-1)]\otimes 
\C[L_C]$
of $V$ is zero.
 
Recall that $V=\sum_{n\geq 0}V_n$ with $V_0=\C[L_C].$ By the theory
of $A(V)$ [Z], $\C[L_C]$ is a simple $A(V)$-module such that
$o(u)=u_{\wt u-1}$ for a homogeneous $u\in V.$ In particular,
$o(e^{\alpha})$ acts on $V_0$ as the multiplication by $e^{\alpha}$
for $\alpha\in L_C$ and $o(d_i(-1))=d_i(0)$ acts on $e^{\beta}$ for
$\beta\in L_C$ as scalar $(d_i,\alpha).$ If we identify 
$\C[L_C]$ with the ring $\C[t_1,t_1^{-1},...,t_{\nu},t_{\nu}^{-1}]$
by identifying $e^{\sum_in_ic_i}$ with $\prod_{i}t_i^{n_i}$ then
$o(e^{\sum_in_ic_i})$ acts on 
$\C[t_1,t_1^{-1},...,t_{\nu},t_{\nu}^{-1}]$
as multiplication by $\prod_it^{n_i}$ and $o(d_i(-1))$ acts 
as the degree derivation $t_i\frac{\partial}{\partial t_i}.$ It is
immediate then that if $o(v)=0$ on $V_0$ for $v\in \C[d_1(-1),..., 
d_{\nu}(-1)]\otimes \C[L_C]$ then $v=0,$ as desired.
\qed
 
We next use $A(V)$ to study the $\Z_+$-graded modules for $V.$ A 
$V$-module
$M=(M,Y_M)$ is $\Z_+$-graded if $M=\oplus_{n\geq 0}M(n)$
such that $u_nM(m)\subset M(\wt u-n-1+m)$ for homogeneous
$u\in V$ and $m,n\in \Z.$ We may and will choose the gradation
so that $M(0)\ne 0.$ Then the theory of $A(V)$ [Z] says 
that $M(0)$ is a $A(V)-module$ by sending $v$ to $o(v).$ Furthermore,
$M\mapsto M(0)$ gives a bijection between the equivalence classes 
of irreducible $\Z_+$-graded $V$-modules and the equivalence
classes of simple $A(V)$-modules. 
 
Recall from Section 4 that given an $A$-module $W$ and $\lambda\in 
\frac{1}{k}L_D$ we have a $V$-module $V_{\l,W}.$ On the other hand, by the
theory in [Z] there is an $\Z_+$-graded $V$-module
$M=\oplus_{n\geq 0}M(n)$ such that $M(0)$ is isomorphic to $W$ as
$A(V)$-modules. If $W$ is simple then $V_{\l,W}$ is irreducible
$V$-module and we can choose $M$ to be irreducible too. Here is the 
relation
between $V_{\l,W}$ and $M$ if $W$ is simple. In the construction
of the module $V_{\l,W}$ we let $c_i(0)$ act on $V_{\l,W}$ as $(c_i,\l)$
($c_i$ is not an element of $A$). On the other hand, since $c_i(-1)\in 
O(V)$
we see that $c_i(0)$ acts on $M=0.$ Thus $V_{0,W}$ and
$M$ are  isomorphic. 
 
Of course we could define the algebra $A$ by including the central 
elements
$c_i.$ Then $A(V)$ would be a quotient of $A$ modulo the ideal
generated by $c_i$ for $i=1,...,\nu.$

\end{document}